\documentclass[11pt]{amsart}
\usepackage[utf8]{inputenc}
\usepackage[pdftex]{graphicx}
\usepackage{amssymb}
\usepackage{amsmath}
\usepackage{amsfonts}
\usepackage{amsthm}
\usepackage{bm}
\usepackage{mathrsfs}
\usepackage{appendix}
\usepackage{enumerate}
\usepackage[left=3cm, right=3cm, bottom=3.4cm]{geometry}
\usepackage[stretch=30,shrink=30]{microtype}
\usepackage{stmaryrd}
\usepackage[normalem]{ulem} 

\usepackage{xcolor}
\definecolor{antiquefuchsia}{rgb}{0.57, 0.36, 0.51}
\definecolor{azure}{rgb}{0.0, 0.5, 1.0}

\usepackage[backref=page, colorlinks = true, urlcolor = black,  linkcolor = black, citecolor = antiquefuchsia]{hyperref}

\renewcommand*{\backref}[1]{}
\renewcommand*{\backrefalt}[4]{%
    \ifcase #1 (Not cited.)%
    \or        (Cited on page~#2.)%
    \else      (Cited on pages~#2.)%
    \fi}

\makeatletter
\def\th@plain{%
	\thm@notefont{}
	\itshape 
}
\def\th@definition{%
	\thm@notefont{}
	\normalfont 
}
\makeatother

\numberwithin{equation}{section}

\newtheorem{theorem}{Theorem}[section]
\newtheorem{lemma}[theorem]{Lemma}
\newtheorem{proposition}[theorem]{Proposition}
\newtheorem{corollary}[theorem]{Corollary}

\theoremstyle{definition}

\theoremstyle{remark}
\newtheorem{remark}[theorem]{Remark}

\newcommand{\N}{\mathbb{N}}

\newcommand{\R}{\mathbb{R}}
\renewcommand{\S}{\mathbb{S}}

\usepackage[capitalise,nameinlink]{cleveref}
\crefformat{equation}{#2(#1)#3}
\crefrangeformat{equation}{#3(#1)#4 to~#5(#2)#6}
\crefmultiformat{equation}{#2(#1)#3}{ and~#2(#1)#3}{, #2(#1)#3}{ and~#2(#1)#3}

\title[]{Quantitative stability of Sobolev inequalities on compact Riemannian manifolds}
\author[Francesco Nobili]{Francesco Nobili} 
\address{Universit\'a di Pisa, Dipartimento di Matematica, Largo Bruno Pontecorvo 5,
56127 Pisa, Italy}
\email{\url{francesco.nobili@dm.unipi.it}}

\author[Davide Parise]{Davide Parise} 
\address{University of California San Diego, Department of Mathematics, 9500 Gilman Drive \#0112, La Jolla, CA 92093-0112, United States of America} 
\email{\url{dparise@ucsd.edu}}

\date{\today. }

\begin{document}
\begin{abstract}
    We study quantitative stability results for different classes of Sobolev inequalities on general compact Riemannian manifolds. We prove that, up to constants depending on the manifold, a function that nearly saturates a critical Sobolev inequality is quantitatively $W^{1,2}$-close to a non-empty set of extremal functions, provided that the corresponding optimal Sobolev constant satisfies a suitable strict bound. The case of sub-critical Sobolev inequalities is also covered. Finally, we discuss degenerate phenomena in our quantitative controls.
\end{abstract}
\maketitle
\tableofcontents

\section{Introduction} 

In the standard round sphere $\mathbb{S}^d$, for $d>2$, the sharp Sobolev inequality reads:
\begin{equation}
\|u\|_{L^{2^*}}^2\le {\sf S}_d^2 \| \nabla u\|_{L^{2}}^2 + { \rm Vol}(\S^d)^{-2/d}\| u\|_{L^2}^2, \qquad \forall u \in W^{1,2}(\S^d),
\label{eq:Sob sphere}
\end{equation}
where $2^* := 2d /(d-2)$ is the Sobolev conjugate exponent, and ${\sf S}_d >0$ is the optimal constant in the sharp Euclidean Sobolev inequality
\begin{equation}
\|u\|_{L^{2^*}}\le {\sf S}_d\|\nabla u\|_{L^2},\qquad \forall u \in \dot W^{1,2}(\R^d).\label{eq:Sob Eucl}
\end{equation}
The value
\[
 {\sf S}_d = \left(\frac{2^*-2}{d}\right)^{\frac 12}{\rm Vol}(\S^d)^{-1/d},
\]
was found by Aubin \cite{Aubin76-2} and Talenti \cite{Talenti76}, while the validity of \eqref{eq:Sob sphere} goes back to \cite{Aubin76-3} where also non-constant \emph{extremal functions} were classified (see \cite[Chapter 5]{Hebey99}):
\begin{equation}
 u_{a,b,z_0} := a(1-b\cos(d(\cdot ,z_0))^{\frac{2-d}{2}},\qquad \text{with } a\in\R,\, b\in(0,1),\, z_0 \in \S^d.
\label{eq:spherical bubbles}
\end{equation}
Here $d$ is the geodesic distance on $\S^d$. In light of this, a natural question is the one of stability:
\begin{center}
If $u\neq cst$ nearly \emph{saturates} \eqref{eq:Sob sphere}, is $u$ close to some  $u_{a,b,z_0}$?
\end{center}
This question is equivalent (up to a change of coordinates via the stereographic projection, see e.g.\ \cite[Appendix A]{FrankLieb12}) to the analogous question raised in \cite{BLi} for the Euclidean Sobolev inequality \eqref{eq:Sob Eucl} and was addressed in Bianchi and Egnell's work \cite{BianchiEgnell91}. The \emph{quantitative stability} for \eqref{eq:Sob Eucl} initiated in \textit{loc. cit.} has been then studied thoroughly in subsequent works \cite{CianchiFuscoMaggiPratelli09,FigalliNeumayer19,Neumayer19,FigalliZhang22} dealing with $p$-Sobolev inequalities for $p\neq 2$, and recently in \cite{DolbeaultEstebanFigalliFrankLoss23,DolbeaultEstebanFigalliFrankLoss24} with explicit constant (see also to the related \cite{Konig23} and  \cite{Frank23_Survey}).  \smallskip

\noindent Aim of this manuscript is to study quantitative stability of optimal Sobolev inequalities on a general closed Riemannian manifold in connection with the existence of extremal functions. We now start by recalling the notions of optimal constants appearing in the celebrated AB-program referring to \cite{Hebey99}, see also \cite{Aubin98,DruetHebey02} and references therein. Given $d>2$ and $(M,g)$ a closed, i.e. compact and boundaryless, $d$-dimensional smooth Riemannian manifolds, we can consider for constants $A,B\ge 0$ the following type of Sobolev inequalities
\begin{equation}
\|u\|_{L^{2^*}}^2 \le A\|\nabla u\|_{L^2}^2 + B \| u\|_{L^2}^2,\qquad \forall u \in W^{1,2}(M).
\label{eq:SobAB}\tag{{\color{red}$\ast$}}
\end{equation}
Notice that an equation as \eqref{eq:Sob Eucl} cannot hold due to the presence of constant functions and, differently from the Euclidean case, it is not straightforward to consider a single sharp Sobolev inequality due to the presence of two constants. The AB-program then starts with the definition of the following Sobolev constants:
\[
\begin{array}{cc}
\alpha(M):= \inf \{ A \ : \  \eqref{eq:SobAB} \text{ holds for some $B$}  \}, &\quad  \beta(M) := \inf \{ B \ : \  \eqref{eq:SobAB} \text{ holds for some $A$} \},
\end{array}
\]
where we understand the infimum over an empty set as being $+\infty$. The first natural problem is then to determine the value of $\alpha(M)$ and $\beta(M)$. It is rather straightforward to see (e.g.\ \cite[Sec 4.1]{Hebey99}) that the latter satisfies
\[
\beta(M) = {\rm Vol}_g(M)^{-2/d}.
\]
On the other hand, the value of $\alpha(M)$ turns out to be linked to the sharp Euclidean Sobolev constant. More precisely, in \cite{Aubin76-2} (see also \cite[Theorem 4.5]{Hebey99}) it is shown that
\begin{equation}\label{eq:alfa manifold}
	\alpha(M)= {\sf S}_{d}^2.
\end{equation}
In particular, $\alpha(M)$ does not depend on $M$ but only on its dimension $d$. A more subtle question is whether these two constants are actually attained, i.e. if they are minima. For instance, for $M=\S^d$ we have in \eqref{eq:Sob sphere} the validity of a Sobolev inequality with $A=\alpha(\S^d)$ and $B = \beta(\S^d)$ and they are attained by the family \eqref{eq:spherical bubbles}. 

More generally, and starting from the easiest one, it was shown in \cite{Bakry94} that $\beta(M)$ is a minimum (see also \cite[Theorem 4.2]{Hebey99}), i.e.\ there is $\overline A>0$ so that \eqref{eq:SobAB} holds with $\overline A$ and $\beta(M)$. It is instead a deep result of \cite{HebeyVaugon96} that $\alpha(M)$ is attained, i.e.\ there is $\overline B>0$ so that \eqref{eq:SobAB}  holds with $\alpha(M)$ and $\overline B$. Thanks to these results, we can define two further notions of Sobolev constants. We can fix $B = \beta(M)$ and proceed with the $A$-part of this program, i.e.\ minimizing over the admissible $A>0$ and consider more generally all sub-critical ranges 
\begin{equation}
\|u\|_{L^{q}}^2 \le A\|\nabla u\|_{L^2}^2 + {\rm Vol}_g(M)^{2/q-1} \| u\|_{L^2}^2,
\label{eq:Sob_qA}\tag{{\color{red}$\ast_A$}}
\end{equation}
for $q \in (2,2^*]$. Conversely, we can fix $A = \alpha(M)$ and proceed with the $B$-part of this program, i.e.\ minimizing over the admissible $B>0$ and consider
\begin{equation}
\|u\|_{L^{2^*}}^2 \le {\sf S}_d^2\|\nabla u\|_{L^2}^2 + B \| u\|_{L^2}^2.
\label{eq:Sob_qB}\tag{{\color{red}$\ast_B$}}
\end{equation}
Then, we can define the corresponding notions of second-best optimal constants
\[
A_q^{\rm opt}(M) := \inf \{A>0 \colon  \eqref{eq:Sob_qA}\text{ holds with }A\},\quad B_{2^*}^{\rm opt}(M)  := \inf \{B>0 \colon  \eqref{eq:Sob_qB}\text{ holds with }B\},
\]
Differently from $\alpha(M),\beta(M)$ these two constants are automatically minima. Moreover, universal bounds depending on the geometry of $M$ are given for $ A_q^{\rm opt}(M)$  in \cite[Theorem 4.4]{DH02} (see also \cite[Proposition 5.1]{NobiliVIolo22} for general $q$), and for $ B_{2^*}^{\rm opt}(M)$ in \cite{HebeyVaugon95}.\\
%
%
Having established two classes of optimal Sobolev inequalities, it is natural to investigate the existence of extremal functions. Let us consider the following sets:
\begin{align*}
    \mathcal{M}_q(A) &:= \{u \in W^{1,2}\setminus\{0\}  \colon  \text{equality holds with } u \text{ in \eqref{eq:Sob_qA} for }A=A^{\rm opt}_q(M) \},\\
    \mathcal{M}_{2^*}(B)& := \{u \in W^{1,2}\setminus\{0\}  \colon   \text{equality holds with } u \text{ in \eqref{eq:Sob_qB} for }B=B^{\rm opt}_{2^*}(M) \}.
\end{align*}
The existence and compactness properties of extremal functions become challenging in the critical case due to the apriori lack of compactness of the embedding $W^{1,2}\hookrightarrow L^{2^*}$. We will see that compactness properties are related to the values of the constants $A^{\rm opt}_{2^*}(M),B^{\rm opt}_{2^*}(M)$. Notice however that $\mathcal M_q(A)$ always contains constant functions thanks to the value of the $B$-constant in \eqref{eq:Sob_qA}. In particular, it is never empty as defined. However, we do not restrict our analysis to \emph{non-constant} extremal functions in this case, as these could be the only extremal functions (see \cite{Frank21} and below). Finally, we refer the reader to  \cite{Hebey01,DjadliDruet01,BarbosaMontenegro07,BarbosaMontenegro09,Barbosa10,BarbosaMontenegro12}, where extremal functions in the AB-program have been investigated. 

\smallskip
\noindent\textbf{Statement of the main results}. 
Here we present our main quantitative stability results. As already observed, $\mathcal M_q(A)$ is never empty as it contains constant functions. Also, by definition of $\alpha(M)$, we always have $A^{\rm opt}_{2^*}(M) \ge {\sf S}_d^2$. However, assuming this inequality to be strict, we gain pre-compactness of (normalized) extremizing sequences via concentration-compactness methods in the spirit of \cite{Lions84,Lions85}, cf. Proposition \ref{prop:precompactness extremal}. Under this assumption, or more easily in the sub-critical range, we can prove the following. 
\begin{theorem}\label{thm: stability Aopt intro}
    Let $(M,g)$ be a closed $d$-dimensional Riemannian manifolds, $d>2$. Assume that either $2<q<2^*$, or $q=2^*$ and $A^{\rm opt}_{2^*}(M)> {\sf S}_d^2$. Then, there are non-negative constants $C,\gamma$ depending on $(M,g)$ and on $q$ so that, for every $u \in W^{1,2}(M)\setminus\{0\}$, it holds
    \[
        \frac{ A^{\rm opt}_q(M)\|\nabla u\|^2_{L^2} +{\rm Vol}_g(M)^{2/q- 1}\|u\|_{L^2}^2}{\| u\|_{L^{q}}^2} - 1 \ge C \left( \inf_{v \in {\mathcal M}_q(A)} \frac{\| u- v\|_{W^{1,2}} }{\|u\|_{W^{1,2}}} \right) ^{2+\gamma}.
    \]
\end{theorem}
The above stability result in Theorem \ref{thm: stability Aopt intro} is \emph{strong}, as we are able to control the full $W^{1,2}$-distance from a set a non empty set of extremal functions with the related Sobolev deficit. On the other hand, the constant $C$ and the exponent $\gamma$ crucially depend on the manifold $M$ and on the value $q$ as an outcome of the proof-strategy. A natural question at this point is whether, on a manifold $(M,g)$, we can expect the strict binding inequality to hold and, therefore, a critical stability result. This heavily depends on the metric $g$, as it is known  from \cite[Corollary 2]{DjadliDruet01} that we always have $A^{\rm opt}_{2^*}(M,h)={\sf S}_d^2$ (adding the dependence on the metric) for a suitable metric $h$ in the conformal class of $g$. However, for certain manifolds this scenario is very rigid (see Proposition 5.9 and Proposition 5.12 in \cite{Hebey99}). More generally, the continuity properties studied in \cite[Theorem 6.2]{NobiliVIolo22} suggests that, as soon as there is $g$ satisfying $A^{\rm opt}_{2^*}(M,g)>{\sf S}_d^2$, there are infinitely many others metrics.

We pass now to our second main result, which is a weaker stability results only around minimizer for $B^{\rm opt}_{2^*}(M)$. Let us first comment on the non-emptyness of the set of extremal functions for $B^{\rm opt}_{2^*}(M)$. When $d\ge4$, work of Hebey \cite[Proposition 5.1]{Hebey99} implies 
\begin{equation} \label{eq:B lower bound}
    B^{\rm opt}_{2^*}(M) \ge \frac{d-2}{4(d-1)}{\sf S}_d^2\max_M {\rm R}_g,
\end{equation}
where ${\rm R}_g$ denotes the scalar curvature of $M$. When this inequality is strict, the deep analysis of \cite{DjadliDruet01} guarantees that $\mathcal M_{2^*}(B)\neq \varnothing$. See also Theorem \ref{thm:dichotomy for Bopt} below for further details. As constant functions might not be extremals, unless $B^{\rm opt}_{2^*}(M)=\beta(M)$, this result is striking. The dimensional restriction $d\ge4 $ is due to the need to consider a more general class of $p$-Sobolev inequalities with $p\neq 2$ when $d=3$, see \cite{DjadliDruet01}. Under this assumption, we prove the following. 
\begin{theorem}\label{thm: stability Bopt intro}
      Let $(M,g)$ be a closed $d$-dimensional Riemannian manifolds, $d\ge 4$. Assume $B^{\rm opt}_{2^*}(M) >\frac{d-2}{4(d-1)}{\sf S}_d^2\max_M {\rm R}_g$. Then, there are non-negative constants $C,\gamma,\delta$ depending on $(M,g)$ so that it holds
    \[
       {\sf S}_d^2\|\nabla u\|^2_{L^2} +B^{\rm opt}_{2^*}(M)\|u\|_{L^2}^2 - 1 \ge C \left( \inf_{v \in {\mathcal M}_{2^*}(B)} \frac{\| u- v\|_{W^{1,2}} }{\|u\|_{W^{1,2}}} \right) ^{2+\gamma},
    \]
    for all $u \in W^{1,2}(M)$ so that $\|u\|_{L^{2^*}} =1$ and $\| u- v\|_{W^{1,2}}  < \delta \|u\|_{W^{1,2}}$ for all $v \in  {\mathcal M}_{2^*}(B)$ with $\|v\|_{L^{2^*}}=1$. 
\end{theorem}
The above is a weaker form of stability due to the assumption forcing $u$ to be close to the set of extremizers. Indeed, contrary to Theorem \ref{thm: stability Aopt intro}, here we cannot rely on a concentration-compactness principle (see Remark \ref{rmk: no CC for Bopt} for details). Besides, the continuous dependence $g\mapsto B^{\rm opt}_{2^*}(M,g)$ studied in \cite{BarbosaMontenegro12} implies the abundance of metrics for which \eqref{eq:B lower bound} holds strict.

\medskip
\noindent Given the above results, it is natural to ask if our results hold \emph{sharp} with quadratic exponent, i.e.\ with $\gamma=0$. This is the typical desirable feature in stability problems. For instance, in \cite{ChodoshEngelsteinSpolaorSpolaor23}, respectively \cite{EngelsteinNeumayerSpolaor22}, the authors show that quadratic stability might fail for the quantitative isoperimetric inequality, respectively the quantitative stability of minimizing Yamabe metrics, on specific manifolds. This degenerate phenomenon has been later analyzed by \cite{Frank21} (see Corollary \ref{cor:Frank example} below), and subsequently in \cite{FrankPeteranderl24}, and \cite{BrigatiDolbeaultSimonov23,BrigatiDolbeaultSimonov24,BrigatiDolbeaultSimonov24_2} for different functional inequalities. Our next results deal with degenerate phenomena.
\begin{theorem}\label{thm:degenerate Aopt main}
    Let $q \in (2,2^*]$ and let $(M,g)$ be satisfying the hypothesis of Theorem \ref{thm: stability Aopt intro}. Assume further that there are no non-constant extremal functions, i.e.\ $\mathcal M_q(A) = \{ c \colon c\in \R\setminus\{0\}\}.$ Then, the stability in Theorem \ref{thm: stability Aopt intro} is degenerate, i.e.\ it must hold with some $\gamma >0$. 
\end{theorem}
We point out that this result is in line with the examples studied by Frank in \cite{Frank21}, as we prove in the next corollary. The exponent $\gamma=2$ follows from the analysis in \textit{loc. cit.} 
\begin{corollary}\label{cor:Frank example}
    For $d>2$, we have: 
    \begin{itemize}
        \item[{\rm i)}] Theorem \ref{thm: stability Aopt intro} does not hold with $\gamma=0$ for $q \in (2,2^*)$ and $M=\S^d$;
        \item[{\rm ii)}]  Theorem  \ref{thm: stability Aopt intro} does not hold with $\gamma=0$ for $q =2^*$ and $M=\mathbb S^1\left(\tfrac{1}{\sqrt{d-2}}\right)\times \mathbb S^{d-1}$.
    \end{itemize}
Finally, in all of the above, $\gamma=2$ is sharp.
\end{corollary}
It is unclear whether this degenerate phenomenon can happen also in Theorem \ref{thm: stability Bopt intro}, and if $M_*=\mathbb S^1(1/\sqrt{d-2})\times \mathbb S^{d-1}$ provides again a counterexample. Moreover, we are only aware of estimates for $B^{\rm opt}_{2^*}(M_*)$ (\cite{HebeyVaugon92} and \cite[Proposition 5.4]{Hebey99}) and not its explicit value.

\smallskip
\noindent\textbf{Comparison with related works}. Quantitative stability results are ubiquitous in the literature. In the classical Euclidean setting, several geometric and functional inequalities have been analyzed from the viewpoint of stability with different techniques. We refer to the aforementioned references regarding the Sobolev inequality, to the surveys \cite{Figalli_ECM, Fusco15Survey}, and the references therein.

On the contrary, in the non-linear framework of Riemannian manifolds, we cannot rely on underlying symmetries of any kind as, for instance, in \eqref{eq:spherical bubbles} for $\S^d$. In fact, our main argument follows a well-established approach, that is to argue via a Lyapunov-Schmidt reduction argument \cite{Simon83} and \L ojasievicz inequality \cite{Lojasiewicz65} to produce quantitative stability estimates. This fruitful interaction has been pionereed in \cite{ChodoshEngelsteinSpolaorSpolaor23} for the quantitative stability of the Riemannian isoperimetric inequality, and it was later employed to study stability properties of minimizing Yamabe metrics in \cite{EngelsteinNeumayerSpolaor22} (see also the recent \cite{ChenKim24}, we refer to \cite{Yamabe60} and the surveys \cite{LP87,BrendleMarques11} regarding the Yamabe problem). In the present work, we exploit the similarities between the latter and our setting.

Besides degenerate phenomena, it is natural to investigate sharp stability results with exponent $\gamma=0$ holding \emph{generically} in the space of Riemannian metrics. However, in contrast with \cite{EngelsteinNeumayerSpolaor22} where this follows by work Schoen \cite{Schoen89}, see also \cite{Anderson05}, the picture here seems to be more delicate. First, to formulate a generic statement, one needs to check that the strict inequalities on our optimal constants are open conditions in the set of $C^3$-metrics. In this direction, recall the continuity properties of $g\mapsto A^{\rm opt}_q(M,g),g\mapsto B^{\rm opt}_{2^*}(M,g)$ studied in \cite{NobiliVIolo22} and \cite{BarbosaMontenegro12}, respectively. Then, the key point would be to prove that generically inside the relevant open sets of metrics, minimizers (or, more generally, non-negative critical points) of
\[
 \frac{ A^{\rm opt}_q(M,g)\|\nabla u\|^2_{L^2} +{\rm Vol}_g(M)^{2/q- 1}\|u\|_{L^2}^2}{\| u\|_{L^{q}}^2} \quad \text{and} \quad \frac{{\sf S}_d^2\|\nabla u\|^2_{L^2} +B^{\rm opt}_{2^*}(M,g)\|u\|_{L^2}^2}{\| u\|_{L^{2^*}}^2}
\]
are non-degenerate. Notice again that the first always admits constant minimizers while, the second does not in general, unless $B^{\rm opt}_{2^*}(M,g)={\rm Vol}_g(M)^{2/2^*- 1}$ (this, however, can always happen in the conformal class of a given metric again by \cite[Corollary 2]{DjadliDruet01}). Constant critical points are typically excluded from this kind of analysis, see e.g. \cite{GhimentiMicheletti,MichelettiPistoia}. Recall also that they can be pathological in this regard, as discussed in Theorem \ref{thm:degenerate Aopt main} and Corollary \ref{cor:Frank example}. Finally, even restricting our investigation to nonconstant extremal functions, it seems that standard transversality arguments \cite{Henry, White} would require extra regularity on the dependence $g\mapsto A^{\rm}_{q}(M,g),g\mapsto B^{\rm opt}_{2^*}(M,g)$.

\smallskip
\noindent We conclude by mentioning that assuming Ricci curvature constraints makes it possible to investigate more sophisticated stability statements in comparison geometry. In this fashion, we recall \cite{CavallettiMaggiMondino19} for the L\'evy-Gromov isoperimetric inequality, and \cite{BerardBessonGallot85,Bertand07,CavallettiMondinoSemola23,FathiGentilSerres22} for the Lichnerowicz spectral gap. Concerning Sobolev inequalities, \emph{qualitative} stability results were deduced in \cite{NobiliVIolo22,NobiliViolo22-2} under Ricci lower bounds. The main difference with this note (besides the qualitative/quantitative analysis and curvature bounds) is that in \cite{NobiliVIolo22,NobiliViolo22-2} the stability properties are studied facing lack of compactness and with explicit classes of -a priori non-extremal- functions having the radial profile of bubbles. Besides, quantitative stability results in the Hyperbolic space have been obtained for the isoperimetric problem in \cite{BolgeleinDuzaarScheven15}, and for the sharp Poincar\'e-Sobolev inequality in \cite{BhaktaGangulyKarmakarMazumdar22,BhaktaGangulyKarmakarMazumdar23PartI}.
\section{Properties of Sobolev inequalities}
\subsection{Optimal constants and extremal functions}
In this part, we prove a key pre-compactness result for normalized extremizing sequences of $A^{\rm opt}_{2^*}(M)$, and mention analogous compactness properties for extremal functions of $B^{\rm opt}_{2^*}(M)$. We start with the former. 
\begin{proposition}\label{prop:precompactness extremal}
    Let $(M,g)$ be a compact $d$-dimensional Riemannian manifolds, $d>2$. Assume that either $q<2^*$ or 
        \[ 
        q=2^*\qquad \text{and}\qquad A^{\rm opt}_{2^*}(M)> {\sf S}_d^2.
        \]
        Then, for every $(u_n)\subseteq W^{1,2}(M)$ non-constant with $\| u_n\|_{L^q}\equiv 1$ satisfying
        \[
       \frac{\| u_n\|_{L^q}-{\rm Vol}_g(M)^{2/q-1}\|u_n\|_{L^2}}{\|\nabla u_n\|_{L^2}}\to A^{\rm opt}_q(M),\qquad\text{as }n\uparrow\infty,
        \]
        it holds up to subsequence that $u_n$ converges strongly in $L^{q}$ and in $W^{1,2}$ to some $u \in \mathcal M_q(A)$ with $\|u\|_{L^{2^*}}=1$. In particular, $\{ u \in \mathcal M_q(A) \colon \|u\|_{L^q}=1\}$ is pre-compact in the $W^{1,2}$-topology.
\end{proposition}
\begin{proof}
    Up to a not relabelled subsequence, we can assume by Rellich pre-compactness that $u_n \to u$ in $L^q$ for every $q<2^*$. Moreover, by assumption, we have for some $\delta_n\downarrow 0$ that
    \begin{equation}
    1 = \| u_n\|_{L^q}^2 \ge (A^{\rm opt}_q(M)   -\delta_n)\|\nabla u_n\|_{L^2}^2 +  {\rm Vol}_g(M)^{2/q-1} \|u_n\|_{L^2}^2,   
    \label{eq:extremizing A}
    \end{equation}
    In particular, $u_n$ is $W^{1,2}$-bounded, and, up to a further subsequence, we also have $u_n \rightharpoonup u$ weakly in $W^{1,2}(M)$ to some $u \in W^{1,2}(M)$.
    
    If $q<2^*$, by lower-semicontinuity of the gradient term and $L^2$-convergence of $u_n$ to $u$, from \eqref{eq:extremizing A} we reach
    \[
    1 \ge A^{\rm opt}_q(M)  \| \nabla u\|_{L^2} +{\rm Vol}_g(M)^{2/q-1} \|u\|_{L^2} \ge \|u\|_{L^q}^2.
    \]
    However, $ \|u\|_{L^q} =\lim_n\|u_n\|_{L^q}=1$ so that $u \in \mathcal M_q(A)$, giving in turn that $\|\nabla u_n\|_{L^2} \to \|\nabla u\|_{L^2}.$ In particular, there is convergence of $\| u_n\|_{W^{1,2}}$ to $\|u\|_{W^{1,2}}$. Being $u_n$ already weakly converging, this is equivalent to strong $W^{1,2}$-convergence concluding the proof in this case.

    If $q=2^*$ we argue by concentration-compactness. We shall use Proposition \ref{prop:CC} below and make use of the assumption in this case to rule out concentration phenomenona. If (a) in Proposition \ref{prop:CC} occurs, then an analogous argument as done after \eqref{eq:extremizing A} gives the conclusion. We are left to exclude (b) in Proposition \ref{prop:CC} using the assumption in this case. By definition of $\alpha(M)= {\sf S}_d^2$ (cf. \eqref{eq:alfa manifold}), for every $\epsilon>0$ there is $B_\epsilon>0$ so that, for every $n \in\N$ it is possible to write
    \[
    1= \|u_n\|_{L^{2^*}}^2 \le ( {\sf S}_d^2+\epsilon)\| \nabla u_n\|_{L^2}^2 + B_\epsilon \|u_n\|_{L^2}.
    \]
    This and the fact that $\|u_n\|_{L^2}\to 0$ guarantee that $\liminf_n \| \nabla u_n\|_{L^2}>0$. Rearranging now the above with \eqref{eq:extremizing A}, taking  $n$ to infinity, and recalling that $u_n$ is $L^2$-converging to zero in this case, we get
    \[
    A^{\rm opt}_{2^*}(M) \le  {\sf S}_d^2 + \epsilon.
    \]
    We see now that for $\epsilon$ small enough we eventually reach a contradiction with $ A^{\rm opt}_{2^*}(M)> {\sf S}_d^2$. So (b) does not occur and the proof is concluded also in this case.    

    Finally, the last statement follows easily from the above analysis. Indeed, a sequence $(u_n)\subseteq \{ u \in \mathcal M_q(A) \colon \|u\|_{L^q}=1\}$ is either eventually constant, or not. In the first situation, the renormalization guarantees that, up to a non-relabelled subsequence, $u_n \equiv cst.$ and we are done. The second situation is instead a particular case of what was addressed before. The proof is therefore concluded.
\end{proof}
In what follows we prove a dichotomy result based on concentration compactness arguments that were needed in the above proof. This was studied in \cite[Theorem 6.1]{NobiliVIolo22} in the spirit of \cite{Lions84,Lions85} but in a different setting (converging metric measure spaces with synthetic Ricci bounds). We include here a quick proof with the core of the argument to be self-contained.
\begin{proposition}\label{prop:CC}
Let $(M,g)$ be a compact and complete $d$-dimensional Riemannian manifold, $d>2$. Suppose that it holds
\begin{equation}
\|u\|_{L^{2^*}}^2\le A\|\nabla u\|_{L^2}^2+B\|u\|_{L^2}^2,\qquad \forall u \in W^{1,2}(M),
\label{eq:Sob prop AB}
\end{equation}
for some constants $A,B>0$, and there exist a sequence $(u_n) \subset W^{1,2}(M)\setminus \{0\}$ satisfying 
	\[ \|u_n\|^2_{L^{2^*}} \ge A_n\| \nabla u_n\|^2_{L^2} + B_n\|  u_n \|^2_{L^2},\]
	for some sequences $A_n \to A$, $B_n\to B$. Then, setting $\tilde u_n:= u_n\|u_n\|_{L^{2^*}}^{-1}$, there exists a non-relabeled subsequence such that only one of the following holds:
	\begin{itemize}
		\item[{\rm (a)}] $\tilde u_{n}$ converges $L^{2^*}$-strong to a function $u \in W^{1,2}(M)$;
		\item[{\rm (b)}] $\|\tilde u_{n}\|_{L^2}\to 0$ and there exists $p \in M$ so that $|\tilde u_n|^{2^*}{\rm Vol}_g  \rightharpoonup \delta_p$ in duality with $C(M)$.
	\end{itemize}
\end{proposition}
\begin{proof}
  By a scaling argument, we can assume $\|u_n\|_{L^{2^*}} \equiv 1$. As $(A_n,B_n)\to (A,B)$, we directly deduce that $ u_n$ is $W^{1,2}$-bounded. Hence, up to a not relabeled subsequence, we have that $u_n$ converges to some $u \in W^{1,2}(M)$ weakly in $W^{1,2}$ and in $L^{2^*}$, and strongly in $L^2$, to some function $u \in W^{1,2}(M)$ and, by tightness as $M$ is compact, that $|u_n|^{2^*}{\rm Vol}_g\rightharpoonup \nu,$ as well as $|\nabla u_n|^2{\rm Vol}_g \rightharpoonup \mu$, in duality with $C(M)$ for a probability measure $\mu$ and for a finite non-negative measure $\nu$. We have \cite{Lions84,Lions85} (see also the arguments \cite[Lemma 6.6]{NobiliVIolo22}) that there is a countable set of indices $J$, points $(x_j)_{j \in J}\subseteq M$ and weights $(\mu_j)_{j\in J},(\nu_j)_{j \in J}\subseteq \mathbb R^+$, satisfying
    \begin{align*}
        &\nu = |u|^{2^*}{\rm Vol}_g + \sum_{j \in J}\nu_j \delta_{x_j}, \\
        &\mu \ge |\nabla u|^2{\rm Vol}_g + \sum_{j \in J}\mu_j\delta_{x_j},\\
        & \nu_j^{2/2^*} \le A \mu_j,\qquad \sum_{j \in J}\nu_j^{2/2^*} <\infty.
    \end{align*}
    We can therefore estimate
    \begin{align*}
    1 = \lim_{n\uparrow\infty} \|u_n\|^2_{L^{2^*}} &\ge \liminf_{n\uparrow\infty} A_n\| \nabla u_n\|^2_{L^2} + B_n\|  u_n \|^2_{L^2} \\
    &\ge A\mu(M) + B \| u\|_{L^2}^2 \\
    &\ge A\|\nabla u\|_{L^2}^2 + \sum_{j \in J} A\mu_j + B \|u\|_{L^2}^2 \\
    &\ge \Big(\int |u|^{2^*}d{\rm Vol}_g\Big)^{2/2^*} + \sum_{j \in J} \nu_j^{2/2^*} \\
    &\ge \Big( \int |u|^{2^*}d{\rm Vol}_g + \sum_{j \in J}\nu_j \Big)^{2/2^*} = \nu(M)^{2/2^*} = 1,
    \end{align*}
    having used that $u$ satisfies \eqref{eq:Sob prop AB} with $A,B$, the properties of $\mu,\nu$, as well as the fact that $t\mapsto t^{2/2^*}$ is strictly concave. However, since all inequalities must be equalities, we have that either $\|u\|_{L^{2^*}}=1$ and all $\nu_j= 0$ for every $j \in J$, or $u=0$ and there is $j_0 \in J $ so that $\nu_{j_0}=1$ and $\nu_j =0$ for all $j\neq j_0$. The first case is precisely conclusion (a), as weak convergence in $L^{2^*}$ together with convergence of the $L^{2^*}$ norms implies strong convergence. If instead $u=0$, then (b) occurred implying at the same time that $u_n$ is $L^2$-converging to zero, and $|u_n|^{2^*} {\rm Vol}_g$ is weakly converging to a Dirac mass centred at $p = x_{j_0}$. Having analysed all the possibilities, the proof is now concluded.
    \end{proof}
We now discuss the case of $B^{\rm opt}_{2^*}(M)$. Recall that thanks to \cite[Proposition 5.1]{Hebey99}, we have the lower bound \eqref{eq:B lower bound}. Notice that, even arguing with Proposition \ref{prop:CC} for minimizing sequences of $B^{\rm opt}_{2^*}(M)$, we cannot exclude the outcome (b), as this time $A={\sf S}_d^2$ and the convergence to zero of the $L^2$-norm makes it impossible to exploit \eqref{eq:B lower bound}. However, sufficient criteria guaranteeing that $\mathcal{M}_{2^*}(B)$ is not empty have been analysed in \cite[Theorem 1 and Corollary 1]{DjadliDruet01}.
\begin{theorem}\label{thm:dichotomy for Bopt}
    Let $(M,g)$ be a $d$-dimensional compact Riemannian manifold with $d\ge 4$. Then, at least one of the following is true
\begin{itemize}
    \item[{\rm i)}]$\mathcal{M}_{2^*}(B)\neq \varnothing$;
    \item[{\rm ii)}] $B^{\rm opt}_{2^*}(M) = \frac{d-2}{4(d-1)}{\sf S}_d^2\max_M {\rm R}_g.$
\end{itemize}
In particular, if strict inequality holds in \eqref{eq:B lower bound}, then $\mathcal M_{2^*}(B)\neq \varnothing.$ 
\end{theorem}
In lower dimensions, a related result was studied in \cite{DjadliDruet01} by considering $p$-Sobolev inequalities with $p\neq 2$. As we stick to $p=2$ in this note, we do not pursue this any further. The above is a dichotomy result achieved via a delicate study of blow-up phenomena of concentrating sequences of solutions of partial differential equations. In particular, strict inequality in \eqref{eq:B lower bound} guarantees compactness. More generally, we recall that $\{u \in \mathcal M_{2^*}(B)\colon \|u\|_{L^{2^*}}=1\}$ is pre-compact in $L^{2^*}$-topology if strict inequality holds in \eqref{eq:B lower bound}. This result is commonly attributed to the analysis contained in \cite{DjadliDruet01}, see the introduction in \cite{BarbosaMontenegro12}, as well as \cite[Corollary 1.4]{BarbosaMontenegro12} for generalizations that, as a by-product, imply this fact. Therefore, arguing exactly as done in Proposition \ref{prop:precompactness extremal}, we deduce
\begin{equation} \label{eq:preocompact MB}
    \{u \in \mathcal M_{2^*}(B) \colon \|u\|_{L^{2^*}}=1\} \text{ is $W^{1,2}$ pre-compact, if $d\ge 4$ and \eqref{eq:B lower bound} holds strict}. 
\end{equation}
Let us discuss a relevant class of manifolds when the above discussion applies. Notice that i) must  occur in Theorem \ref{thm:dichotomy for Bopt} if either ${\rm R}_g\le 0$, or ${\rm R}_g$ is constant. The first case being clear, we discuss the second one. The case where $(M,g)$ is conformally equivalent to the round sphere has been fully investigated in \cite[Theorem 5.7]{Hebey99}. In the remaining case, the resolution of the Yamabe problem \cite{Yamabe60} by combination of \cite{Trudinger68, Aubin76-3, Schoen84} guarantees that
\[
Y(M)<{\sf S}_d^{-2}, 
\]
having denoted $Y(M)$ the Yamabe constant of $(M,g)$ that is attained by some $u \in W^{1,2}(M)$ 
\[
\|u\|_{L^{2^*}}^2Y(M) = \|\nabla u\|_{L^2}^2 + \frac{d-2}{4(d-1)}{\rm R_g}\|u\|_{L^2},
\]
having used that ${\rm R}_g$ is constant. Thus, there holds
\begin{align*}
     \|u\|_{L^2}^2\frac{B_{2^*}^{\rm opt}(M)}{{\sf S}_d^2} +\|\nabla u\|_{L^2}^2 &\ge \|u\|_{L^{2^*}}{\sf S}_d^{-2} > \|u\|_{L^{2^*}}Y(M) =    \|\nabla u\|_{L^2}^2 + \frac{d-2}{4(d-1)}{\rm R_g}\|u\|_{L^2}. 
\end{align*}
Rearranging, we get that ii) cannot occur.
\subsection{First and second-order variations}
In this part, we study functional properties of abstract Sobolev quotients. This investigation is analogous to \cite[Section 2]{EngelsteinNeumayerSpolaor22}. We start setting up some notation. We fix $q \in (2,2^*]$. We omit the dependence on it of various quantities to lighten the notation. For $\delta>0$, set
\[
\mathcal B := \{ u \in W^{1,2}(M) \colon u\ge 0, \| u\|_{L^q} =1 \},\qquad \mathcal B_\delta(v) := \{ u \in \mathcal B \colon \|u-v\|_{W^{1,2}}\le \delta\}.
\]
Arguing as in \cite[Lemma 2.1]{EngelsteinNeumayerSpolaor22}, $\mathcal{B} \subseteq W^{1,2}$ is a Banach submanifold and the tangent space at $u \in \mathcal B$ is given by 
\[
T_u\mathcal B := \Big\{  v \in W^{1,2}(M) \colon \int u^{q-1}v\, d {\rm Vol}_g =0\Big\}.
\]
Since $W^{1,2}(M)$ is a Hilbert space, we denote by $\pi_{T_u\mathcal B}$ the $L^2$-projection operator defined by
\[
\pi_{T_u\mathcal B}(\varphi) := \varphi - \Big(\int u^{q-1} \varphi \, d {\rm Vol}_g \Big) u.
\]

We compute now the first and second variations of the functional arising by the Sobolev inequality (cf. \cite[Lemma 2.1]{EngelsteinNeumayerSpolaor22}).
\begin{lemma} \label{lem: variations Q}
Let $(M,g)$ be a compact $d$-dimensional Riemannian manifold for $d>2$. Set for $u\in W^{1,2}(M)\setminus\{0\}$ and $A,B>0, q \in (2,2^*]$:
\[
 \mathcal Q_q(u) := \frac{ A \| \nabla u\|_{L^2}^2 + B \|u\|_{L^2}^2 }{\| u\|_{L^{q}}^2}.
\]
 For $u \in \mathcal B$, the first and second variation of $\mathcal Q_q$ at $u$ restricted to $T_u\mathcal B$ are denoted by  $\nabla_\mathcal{B} \mathcal{Q}_q(u),\nabla^2_\mathcal{B} \mathcal{Q}_q(u)$ and are respectively given by
\[
\nabla_\mathcal{B} \mathcal{Q}_q(u)[\varphi] = 2 \int A\nabla u \nabla \varphi +B u\varphi \,d {\rm Vol}_g, 
\]
and
\[
\nabla^2_\mathcal{B} \mathcal{Q}_q(u)[\varphi,\eta] = 2  \int A\nabla \varphi\nabla\eta + B\varphi\eta \, d{\rm Vol}_g  -2(q-1)\mathcal Q_q(u) \int u^{q-2}\varphi\eta\, d{\rm Vol}_g, 
\]
for all $\varphi, \eta \in W^{1, 2}(M)$ (omitting the projections $\pi_{T_u \mathcal B}$ inside the integrals).

Moreover, the following properties hold:
\begin{itemize}
    \item[{\rm i)}] $w\mapsto \frac{\nabla^2_\mathcal{B}\mathcal{Q}_q(w)[\varphi, \cdot]}{\|\varphi\|_{C^{2,\alpha}}}$ is continuous from $C^{2,\alpha}\cap \mathcal{B}\to C^{0,\alpha}$ with uniform modulus of continuity over $\varphi \in C^{2,\alpha}$;
    \item[{\rm ii)}] $w\mapsto \frac{\nabla^2_\mathcal{B}\mathcal{Q}_q(w)[\varphi,\eta]}{\|\varphi\|_{W^{1, 2}}\|\eta\|_{W^{1, 2}}}$ is continuous from $\mathcal{B}\to \R$ with uniform modulus of continuity over $\varphi,\eta \in W^{1,2}$.
\end{itemize}
\end{lemma}
\begin{proof}
We fix any $q \in (2,2^*]$ and $u \in \mathcal B$, and consider any $\varphi \in C^\infty(M), \epsilon \in (-1,1)$. Set 
\[
u_\epsilon := (u+\epsilon \varphi)/I_\epsilon,\qquad I_\epsilon := \| u+\epsilon \varphi\|_{L^q}.
\]
We claim that $u_\epsilon$ is eventually for $|\epsilon|$ small, well defined. Indeed, from the estimate $\big| |x+\epsilon y|^q -|x|^q\big|\le q|\epsilon y|\big| |x+\epsilon y|^{q-1} + |x|^{q-1}\big|$ together with the fact that $u,\varphi \in L^q$, we reach that $I_\epsilon \to 1$ as $\epsilon \to 0$. For future use, we notice
    \begin{equation}
       \lim_{\epsilon \to 0}\frac{1-I_\epsilon^2}{\epsilon} = \frac{2}{q} \Big(\int |u|^q\,d {\rm Vol}_g \Big)^{\frac{2}{q}-1} \lim_{\epsilon} \frac{ \int |u|^q - |u+\epsilon \varphi|^q\,d{\rm Vol}_g}{\epsilon} = -2\int u^{q-1}\varphi \, d {\rm Vol}_g,\label{eq:Ieps expansion}
    \end{equation}
having used dominated convergence theorem combined with the same estimate as before for the last equality (recall also that $\|u\|_{L^q}=1$). We start therefore computing the first variation
\[
\nabla \mathcal Q_q(u)[\varphi] := \lim_{\epsilon\to 0}\frac{\mathcal Q_q(u_\epsilon) - \mathcal Q_q(u)}{\epsilon}.
\]
By polarization, we derive
\begin{align*}
    \frac{\mathcal Q_q(u_\epsilon) - \mathcal Q_q(u)}{\epsilon} &= \frac{ \int A |\nabla (u + \epsilon \varphi) |^2 + B |u + \epsilon \varphi|^2\,d {\rm Vol}_g }{\epsilon I^2_{\epsilon}}  - \frac{ \int A|\nabla u|^2 + B |u|^2\, d {\rm Vol}_g}{\epsilon} \\
    &= \frac{1}{\epsilon}\Big[ \frac{1}{I_\epsilon^2} -1\Big] \mathcal Q_q(u) + \frac{2}{I_\epsilon^{2}}\int  A\nabla u\nabla \varphi + Bu\varphi \, d {\rm Vol}_g \\
    &\qquad +  \frac{\epsilon }{I_\epsilon^{2}} \int  A|\nabla \varphi|^2 + B \varphi^2 \, d {\rm Vol}_g .
 \end{align*}
Hence, recalling that $I_\epsilon \to 1$ \eqref{eq:Ieps expansion}, we deduce
\[
\nabla \mathcal{Q}_q(u)[\varphi] = 2 \int A\nabla u\cdot \nabla \varphi + Bu\varphi \, d {\rm Vol}_g - 2\mathcal Q_q(u)\int u^{q-1}\varphi \, d {\rm Vol}_g.
\]
So, the restriction to $T_u \mathcal B$ satisfies
\[
\nabla_\mathcal{B}\mathcal Q_q(u)[\varphi] =  2   \int A\nabla u\cdot\nabla\pi_{T_u\mathcal B}(\varphi)  + Bu\pi_{T_u\mathcal B}(\varphi) \, d {\rm Vol}_g.
\]
We compute now the second variation for 
\[
\nabla^2 \mathcal Q_q(u)[\varphi,\varphi] := \frac{d^2}{d\epsilon^2} \Big|_{\epsilon =0}\mathcal Q_q(u_\epsilon).
\]
Starting from the computations already performed for the first variation, we have
\begin{align*}
\frac{ \mathcal Q_q(u_\epsilon) - 2\mathcal Q_q(u) + \mathcal Q_q(v_{-\epsilon})}{\epsilon^2} &=
  \frac{1}{\epsilon^2}\Big(\frac{1}{I^2_\epsilon} -2 + \frac{1}{I^2_{-\epsilon} }\Big)\mathcal Q_q(u) \\
  &\qquad + \frac{2 }{\epsilon}\Big(\frac{1}{I^2_\epsilon} -\frac{1}{I^2_{-\epsilon}}\Big)\int A\nabla u\nabla \varphi + Bu\varphi \, d {\rm Vol}_g \\
  &\qquad + \Big( \frac{1}{I^2_\epsilon} + \frac{1}{I^2_{-\epsilon}}\Big)\int A |\nabla \varphi|^2 + B \varphi^2\, d {\rm Vol}_g.
\end{align*}
Hence, taking the limit, we reach
\begin{align*}
\nabla^2 \mathcal Q_q(u)[\varphi,\varphi]  
&= -2(q-1)\mathcal Q_q(u)\int u^{q-2}\varphi^2\, d {\rm Vol}_g\\
&\qquad + 8\Big(\int u^{q-1}\varphi\, d {\rm Vol}_g\Big) \int A\nabla u\cdot\nabla \varphi + Bu\varphi \, d {\rm Vol}_g \\
&\qquad + 2\int A|\nabla \varphi|^2 + B \varphi^2\, d {\rm Vol}_g.
\end{align*}
So, the restriction to $T_u \mathcal B$ reads
\begin{align*}
    \nabla^2_\mathcal{B}\mathcal Q_q(u)[\varphi,\varphi]  &= 2  \int A |\nabla \pi_{T_u\mathcal B}(\varphi)|^2 + B \pi_{T_u\mathcal B}(\varphi)^2\, d {\rm Vol}_g \\
    &\qquad-2(q-1)\mathcal Q_q(u)\int u^{q-2}\pi_{T_u\mathcal B}(\varphi)^2\, d {\rm Vol}_g.
\end{align*}
Recalling that $\mathcal{Q}_q(1) = cst.,$ the desired second variation formula hold for all couples $\varphi,\eta \in W^{1,2}(M)$ by polarization.

We now verify i). For every $w \in C^{2,\alpha}$, we define $L_w(\varphi) := -2A\Delta \varphi + B\varphi -2(q-1)\mathcal{Q}_q(w)w^{q-2}\varphi$. Then, for every $w,v \in C^{2,\alpha}$ we have for some $C>0$
\[
\| L_w(\varphi) - L_v(\varphi)\|_{C^{0,\alpha}} \le C \big( \|w^{q-2}\varphi\|_{C^{0,\alpha}}|\mathcal{Q}_q(w)-\mathcal{Q}_q(v)| + \|\varphi |w^{q-2} - v^{q-2}|\|_{C^{0,\alpha}}\big).
\]
Hence, taking into account the continuity of $t\mapsto t^{q-2}$ and that of $w \mapsto \mathcal{Q}_q(w)$ in $C^{2,\alpha}$, it holds
\[
\| L_w(\varphi) - L_v(\varphi)\|_{C^{0,\alpha}} \le f(\| w-v\|_{C^{2,\alpha}})\|\varphi\|_{C^{2,\alpha}},
\]
for some modulus of continuity $t\mapsto f (t)$ with $f(t)\downarrow 0$ as $t\downarrow 0$ uniform on $\varphi$. To prove i), we thus need to prove continuity properties of $w \mapsto  \pi_{T_w \mathcal B}$ as a map from $C^{2,\alpha} \cap \mathcal B\to \mathcal{L}(C^{2,\alpha},C^{2,\alpha})$. This is implied by the estimate
\[
\big| \pi_{T_w\mathcal B}(\varphi) - \pi_{T_v\mathcal B}(\varphi)\big|  = \Big| w-v + \int (v^{q-1}-w^{q-1})\varphi\, d{\rm Vol}_g\Big| \le C \| w-w\|_{C^{0,\alpha}}\|\varphi\|_{C^{2,\alpha}}.
\]
We now prove ii) and conclude. For every $\varphi \in W^{1,2}(M)$, and interpreting by slight abuse of notation $L_w(\varphi)$ in distributional sense, we estimate similarly as before 
\[
\big\| L_w(\varphi) - L_v(\varphi) \big\|_{H^{-1} }\le  C \big( \|w^{q-2}\varphi\|_{H^{-1}}|\mathcal{Q}_q(w)-\mathcal{Q}_q(v)| + \|\varphi |w^{q-2} - v^{q-2}|\|_{H^{-1}}\big).
\]
Then, by H\"older and the Sobolev inequality, we get 
\[
\|w^{q-2}\varphi\|_{H^{-1}} \le \|w^{q-2}\|_{L^{q/(q-2)}} \|\varphi\|_{L^q} \le C \| w\|_{L^q}^{q-2}  \|\varphi\|_{W^{1,2}},
\]
for some $C>0$ and, using $|a-b|^{\frac{1}{q-2}} \le c_q |a^\frac{1}{q-2}- b^\frac{1}{q-2}|$, we estimate similarly
\[
\|\varphi |w^{q-2} - v^{q-2}|\|_{H^{-1}} \le c_q \| w-v\|_{L^q}^{q-2} \|\varphi\|_{L^q} \le f\big( \|w-v\|_{W^{1,2}} \big) \|\varphi\|_{W^{1,2}},
\]
for some modulus of continuity $t\mapsto f (t)$ with $f(t)\downarrow 0$ as $t\downarrow 0$ uniform on $\varphi$. Recall also that $\mathcal Q_q(\cdot)$ is continuous in $W^{1,2}$. From here, the conclusion follows by continuity properties of $w\mapsto \pi_{T_w\mathcal B}$ as a mapping  $\mathcal B \to \mathcal L(W^{1,2},W^{1,2})$. 
\end{proof}
We also have the crucial property that $\mathcal Q_q(\cdot)$ is analytic in the sense of \cite[Lemma 6]{CarlottoChodoshRubinstein2015} whose proof carries over in our setting.
\begin{proposition}\label{prop:analitic quotient}
    Under the same assumptions of Lemma \ref{lem: variations Q}, the map $u \mapsto \mathcal{Q}_q(u)$ is analytic on $\{u \in C^{2, \alpha}(M), u > 0 \}$, in the sense that for every $u_0 \in C^{2, \alpha}(M)$ with $u_0 > 0$, there is $\epsilon > 0$, and bounded multilinear operators $\mathcal{Q}_{q}^{(k)} \colon C^{2, \alpha}(M)^{\times k} \rightarrow \mathbb{R}$ for every $k \geq 0$, such that if $\Vert u - u_0 \Vert_{C^{2, \alpha}(M)} < \epsilon$, then $$\sum_{k = 0}^{\infty} \Vert \mathcal{Q}_{q}^{(k)} \Vert \cdot \Vert u - u_0 \Vert_{C^{2, \alpha}(M)}^{k} < \infty, $$ as well as 
    \begin{equation*}
        \mathcal{Q}_q(u) = \sum_{k = 0}^\infty \mathcal{Q}_{q}^{(k)}(u - u_0, u - u_0, \ldots, u - u_0), \qquad \text{in }C^{2, \alpha}(M).
    \end{equation*}
\end{proposition}
\subsection{Lyapunov-Schmidt reduction} 
The aim of this section is to discuss the Lyapunov-Schmidt reduction associated with a general Sobolev quotient that, intuitively, reduces the infinite dimensional problem to a more manageable finite dimensional one. 

Suppose that
\begin{equation}
    \inf_{u \in W^{1,2}(M)\setminus \{0\}}\mathcal Q_q(u) = 
\inf_{u \in W^{1,2}(M)\setminus \{0\}}\frac{ A \| \nabla u\|_{L^2}^2 + B \|u\|_{L^2}^2 }{\| u\|_{L^{q}}^2} = 1.\label{eq:sharp generic AB}
\end{equation}
for $A,B>0$ and $ q \in (2,2^*]$ either as in Theorem \ref{thm: stability Aopt intro} or Theorem \ref{thm: stability Bopt intro}. We shall not specify their value to unify the discussion here. Consider the set of normalized extremal functions for \eqref{eq:sharp generic AB}
\[
\mathcal M^1 := \{ u \in\mathcal B \colon \mathcal Q_q(u) =1\},
\]
and the set of normalized critical points $\mathcal{CP}^1 \subseteq W^{1,2}(M)$ as the collection of $u \in \mathcal B$ so that 
\[
\int_M A\nabla u\cdot \nabla \varphi + B u\varphi\, d {\rm Vol}_g = \mathcal{Q}_q(u)\int_M u^{q-1}\varphi\, d {\rm Vol}_g, \qquad\forall \varphi \in C^\infty(M).
\]
The sets $\mathcal M^1$ and $\mathcal{CP}^1$ depend on $q,A,B$ and coincide with those discussed in the Introduction. Later, in the main proofs, these constants will be chosen. Furthermore, we can regard $ \nabla^2_\mathcal{B}\mathcal Q_q(u)[\cdot,\cdot]$ as a mapping $T_u\mathcal B\subseteq W^{1,2}(M) \to H^{-1}(M)$. Since the operator associated with this mapping is elliptic on a compact manifold, we know that
\[
K := {\rm Ker} \, \nabla^2_\mathcal{B}\mathcal Q_q(u),
\]
has finite dimension, say dim$(K)= l<\infty$. We then denote by $K^\perp$ its $L^2$-orthogonal complement in $W^{1,2}(M)$. Given the properties of $\mathcal Q_q(\cdot)$, a proof of the following result follows by \cite[Section 3.12]{SimonHarmonic}, as well as \cite[Appendix A]{EngelsteinNeumayerSpolaor22} replacing the use of \cite[Lemma 2.1]{EngelsteinNeumayerSpolaor22} with Lemma \ref{lem: variations Q}. Recall also that any $v \in \mathcal M^1$ is $C^{2,\alpha}$-regular and it s either strictly positive or negative, by elliptic regularity (see, e.g., \cite{Druet00}). We also refer the reader to the instructive finite dimensional version for smooth functions in \cite[Appendix 3.16.1]{SimonHarmonic}.  
\begin{proposition}[Lyapunov-Schmidt reduction]\label{prop:Lyapunov reduction}
    Let $(M, g)$ be a $d$-dimensional compact Riemannian manifold, $d>2$ with a $C^3$ metric. Suppose that \eqref{eq:sharp generic AB} holds with $A,B>0$ and $q\in (2,2^*]$. Fix $v \in \mathcal M^1$. There are an open neighborhood $U$ of $0$, such that $U \subset K$, and a map $F \colon U \rightarrow K^{\perp}$ satisfying $F(0) = 0$, and $\nabla F(0) = 0$ such that the following hold. 
    \begin{itemize}
        \item[{\rm i)}] Let $\mathfrak{q} \colon U \rightarrow \mathbb{R}$ be the analytic function defined by $\mathfrak q(\varphi) = \mathcal Q_q(v + \varphi + F(\varphi))$. Then 
        \[
            \mathcal{L} := \{v + \varphi + F(\varphi); \; \varphi \in U\} \subset \mathcal{B},
        \]
        and 
        \[
            \nabla_\mathcal{B} \mathcal Q_q(v + \varphi + F(\varphi)) = \pi_K \nabla_\mathcal{B} \mathcal Q_q(v  + \varphi + F(\varphi)) = \nabla \mathfrak q(\varphi). 
        \]
        Furthermore, $\varphi \mapsto \mathfrak q (\varphi)$ is real analytic;
        \item[{\rm ii)}] There exists $\delta > 0$, depending on $v$ such that for any $u \in \mathcal{B}_\delta(v)$, we have $\pi_K(u - v) \in U$. Furthermore, if $u \in \mathcal{CP}^1 \cap \mathcal{B}_\delta(v)$, then 
        \[
            u = v + \pi_K(u - v) + F(\pi_K(u - v));
        \]
        \item[\rm iii)] There exists $C > 0$, such that for all $\varphi \in U$ and $\eta \in K$, we have the bound 
        \[
            \Vert \nabla F(\varphi)[\eta] \Vert_{C^{2, \alpha}} \leq C \Vert \eta \Vert_{C^{0, \alpha}}. 
        \]
    \end{itemize}
\end{proposition}
\section{Quantitative stability}
Aim of this section is to prove first a local quantitative stability, and then pass to the proof of our main quantitative stability results.
\subsection{Local quantitative stability}
We shall continue assuming here the validity of \eqref{eq:sharp generic AB}, without specifying the value of $A,B>0$ and $q\in(2,2^*]$. Recall now $K,F$ from Proposition \ref{prop:Lyapunov reduction} and the definitions of $\mathcal M^1$ and $\mathcal{CP}^1$ given above.

Associated to the Lyapunov-Schmidt reduction, we have the following notions (\cite{AdamsSimon1988} see also \cite{CarlottoChodoshRubinstein2015}). For $\delta>0,v \in \mathcal M^1$ given by Proposition \ref{prop:Lyapunov reduction} and $u \in \mathcal B_\delta(v)$, we define the \emph{Lyapunov-Schmidt projection operator}
\[
u_\mathcal{L} := v + \pi_K(u-v) + F(\pi_K(u-v)).
\]
We can thus distinguish between the following notions:
\begin{itemize}
    \item We say that $v$ is \emph{non-degenerate} if ${\rm Ker} \nabla^2_\mathcal{B} \mathcal Q_q(v) = \{0\}$. Hence, $ u_\mathcal{L} = v $ in this case for every  $u \in \mathcal B_\delta(v)$;
    \item We say that $v \in \mathcal{CP}^1$ is \emph{integrable}, provided that for every $\varphi \in {\rm Ker} \nabla^2_\mathcal{B} \mathcal Q_q(v)$, there is a one-parameter family of functions $(v_t)_{t \in (-\delta,\delta)}$ with $v_0 = v$, and $\frac{d v_t}{dt}|_{t=0}=\varphi$, as well as $v_t \in  \mathcal{CP}^1$ for all $t$ sufficiently small. Conversely, we say that $v$ is non \emph{non-integrable} if this property does not hold. 
\end{itemize}
Arguing exactly as in \cite[Lemma 2.4]{EngelsteinNeumayerSpolaor22} in the present setting, we have that if $v$ is integrable, then $\mathfrak{q}$ as given in Proposition \ref{prop:Lyapunov reduction} is constant in a neighborhood of $\{ 0\} \in K$, and when $v$ is an integrable minimizer also
\begin{equation}
\mathcal M^1 \cap \mathcal B_\delta(v) = \mathcal L,
\label{eq:Ball v integrable}
\end{equation}
where $\mathcal L$ is the set given in Proposition \ref{prop:Lyapunov reduction}. We define, for $v \in \mathcal M ^1$ and $\delta>0$, the local distance
\[
d_\delta(u,\mathcal M^1) := \inf   \left\{ \frac{ \|u-\tilde v\|_{W^{1,2}}}{\|u\|_{W^{1,2}}}   \colon \tilde v \in \mathcal M^1 \cap \mathcal B_\delta(v) \right\}.
\]
Recall also the following inequality due to \cite{Lojasiewicz65}. 
\begin{proposition}[\L ojasievicz inequality] \label{prop: loj ineq}
    Let $\mathfrak q \colon \mathbb{R}^l \rightarrow \mathbb{R}$ be a real analytic map and assume $\nabla \mathfrak q(0) = 0$. Then, there exist $\delta > 0$, $C > 0$, and $\gamma > 0$ depending on $\mathfrak q$ and the critical point $0$, such that for all $z \in B(0, \delta)$, it holds
    \begin{equation*}
        \vert \mathfrak q(z) - \mathfrak q(0) \vert \geq C \left( \inf \{ \vert \varphi - \overline{\varphi} \vert; \; \overline{\varphi} \in B_\delta(0), \nabla \mathfrak q(\overline{\varphi}) = 0 \} \right)^{2 + \gamma} .
    \end{equation*}
\end{proposition}
We can prove now the main local stability estimate. 
\begin{proposition}[Local quantitative stability] \label{prop: loc quant stab}
    Let $(M,g)$ be a closed $d$-dimensional Riemannian manifold, $d>2$. Suppose that \eqref{eq:sharp generic AB} holds with $A,B>0$ and $q\in(2,2^*]$. Fix  $v \in \mathcal M^1$. Then, there exist positive constants $C,\gamma,\delta$ depending on $v$ such that
    \[
    \mathcal{Q}_q(u)- 1 \ge C d_\delta(u,\mathcal M ^1)^{2+\gamma},\qquad \forall u \in \mathcal B_\delta(v).
    \]
    Moreover, if $v$ is integrable or non-degenerate, we may take $\gamma =0$.
\end{proposition}
\begin{proof}
   In this proof, $C$ is a general positive constant depending on $v$ (in particular on $(M,g)$ and the value of $q$) whose value might change from line to line.

   Taking $\delta > 0$ in Proposition \ref{prop:Lyapunov reduction} smaller if necessary, we can ensure that for any $\epsilon > 0$, we have $\Vert u_\mathcal{L} - u \Vert_{W^{1, 2}} \leq \epsilon$, as well as $\Vert u_\mathcal{L} - v \Vert_{W^{1, 2}} \leq \epsilon$. In particular, we can write 
   \[
   \mathcal{Q}_q(u) - 1 = (\mathcal{Q}_q(u) - \mathcal{Q}_q(u_\mathcal{L})) + (\mathcal{Q}_q(u_\mathcal{L}) - 1) =: {\rm I} + {\rm II}.
   \]
   We analyse the two term separately. For the first term, write $u = u_{\mathcal{L}} + u^\perp$, and Taylor expand the functional to get
   \begin{align} \label{eq: computation local stab}
   \begin{split}
       {\rm I}  & = \nabla_\mathcal{B} \mathcal{Q}_q(u_{\mathcal{L}})[u^\perp] + \frac{1}{2}\nabla_{\mathcal{B}}^{2} \mathcal{Q}_q(\zeta)[u^\perp, u^\perp] \\
       & = \frac{1}{2}\nabla_{\mathcal{B}}^{2} \mathcal{Q}_q(v)[u^\perp, u^\perp] + o(1) \Vert u^\perp \Vert_{W^{1, 2}}^2 \geq C \lambda_1 \Vert u^\perp \Vert_{W^{1,2}}^2 =C \|u-u_\mathcal{L}\|_{W^{1,2}}^2,
    \end{split} 
   \end{align}
    where $\lambda_1$ is the smallest non-zero eigenvalue of the Hessian of $\mathcal{Q}_q$. Here $\zeta$ is a geodesic point in $\mathcal B$ between $u$ and $u_\mathcal L$ and $o(1)$ is a quantity that goes to zero as $\|u-v\|_{W^{1,2}}$ goes to zero. For the second equality, we used the fact that $\nabla_\mathcal{B} \mathcal{Q}_q(u_{\mathcal{L}})[u^\perp] = 0$, as well as continuity of the Hessian, cf. Lemma \ref{lem: variations Q} property (ii).  The last inequality holds by taking $\delta > 0$ small enough.

    We can now turn to the second term $\mathrm{II}$ specializing the discussion depending on $v$ being non-degenerate, integrable or non-integrable. First, if $v$ is non-degenerate, we have $u_\mathcal{L} = v$, so II $=0$ and \eqref{eq: computation local stab} immediatly yields the conclusion. When $v$ is integrable, we infer $\mathcal{Q}_q(u_{\mathcal{L}}) = \mathfrak q(0) = Q_q(v)$, where the first equality follows from $\mathfrak q$ being constant in a neighborhood of the origin (cf. above \eqref{eq:Ball v integrable}). In particular, also in this situation II $=0$ and since $u_{\mathcal L} \in \mathcal M^1 \cap B_\delta(v)$ the estimate \eqref{eq: computation local stab} combined with the simple observation
    \[
    \|u-u_\mathcal{L}\|_{W^{1,2}}^2 \ge \left(\inf\big\{ \| u - \tilde v\|_{W^{1,2}} \colon \tilde v \in \mathcal M^1\cap B_\delta(v)\big\}\right)^2,
    \]
    immediately give the conclusion.  It remains to handle the case when $v$ is non-integrable. Denoting $\varphi = \pi_K(u - v)$, we deduce
    \begin{align*}
        {\rm II}  = \mathfrak q(\varphi) - \mathfrak q(0) &\geq C \left( \inf \{ \vert \varphi - \overline{\varphi} \vert; \; \overline{\varphi} \in K \cap B_\delta(0), \nabla \mathfrak q(\overline{\varphi}) = 0 \} \right)^{2 + \gamma} \\
        & \geq C \inf \{\Vert u_{\mathcal{L}} - \tilde{v} \Vert_{W^{1, 2}}; \; \tilde{v} \in \mathcal{M}^1 \cap \mathcal B_\delta(v) \}^{2 + \gamma},  
    \end{align*}
    thanks to Proposition \ref{prop: loj ineq}, and where we used the second and third properties of Proposition \ref{prop:Lyapunov reduction}, while for the last inequality we argued as in \cite[page 405]{EngelsteinNeumayerSpolaor22}. Thus, the conclusion of the proof follows by noticing
    \begin{align*}
        \mathcal{Q}_q(u) - 1 &\ge C \|u-u_\mathcal{L}\|_{W^{1,2}}^2 + C \inf \{\Vert u_{\mathcal{L}} - \hat{v} \Vert_{W^{1, 2}}; \; \hat{v} \in \mathcal{M}^1 \cap \mathcal B_\delta(v) \}^{2 + \gamma}\\
        &\ge C\left( \|u-u_\mathcal{L}\|_{W^{1,2}} + \inf \{\Vert u_{\mathcal{L}} - \hat{v} \Vert_{W^{1, 2}}; \; \hat{v} \in \mathcal{M}^1 \cap \mathcal B_\delta(v) \}\right)^{2+\gamma}\\
        &\ge C \left(\inf\big\{ \| u - \tilde v\|_{W^{1,2}} \colon \tilde v \in \mathcal M^1\cap B_\delta(v)\big\}\right)^{2+\gamma}.
    \end{align*} 
\end{proof}
\subsection{Proof of Theorem \ref{thm: stability Aopt intro} and Theorem \ref{thm: stability Bopt intro}}
We start with the proof of Theorem \ref{thm: stability Aopt intro} and then explain how to modify it to obtain Theorem \ref{thm: stability Bopt intro}. 

Fix any $q \in (2,2^*]$ and set here  $A=A^{\rm opt}_q(M)$ and $B={\rm Vol}_g(M)^{2/q-1}$. We claim  that it is enough to prove Theorem \ref{thm: stability Aopt intro} for every $u \in \mathcal B=\{ u \in W^{1,2}(M)\colon \|u\|_{L^q}=1\}$. Indeed, the main estimate is zero-homogeneous in $u$ and 
\[
\inf \{ \| u-v\|_{W^{1,2}} \colon v \in \mathcal M_q(A)\cap \mathcal B\} \ge \inf\left\{ \frac{\| u-v\|_{W^{1,2}}}{\|u\|_{W^{1,2}}} \colon v \in \mathcal M_q(A)\right\} =: d(u,\mathcal M_q(A)).
\]
Thus, fix $v \in \mathcal{M}_q(A)$ with $\Vert v \Vert_{L^q(M)} = 1$.  Let us  invoke Proposition \ref{prop: loc quant stab} with $\mathcal M^1 = \mathcal M_q(A)\cap \mathcal B$ to get the existence of constants $\gamma(v), \delta(v)$, and $C(v)$. In particular, by compactness of $\mathcal{M}_q(A) \cap \mathcal{B}$, c.f. Proposition \ref{prop:precompactness extremal}, we can take a finite subcover (indexed by $I$) of the open cover of balls $\mathcal{B}(v, \delta(v)/2)$. Notice that, in the case $q=2^*$, our extra assumption on $A^{\rm opt}_{2^*}(M)$ plays a crucial role. Set now 
\begin{equation*}
    \delta_0 = \min_{I} \delta(v_i)
/2, \quad \gamma_0 = \max_{I} \gamma(v_i)/2, \quad \text{and} \quad C_0 = \min_{I} C(v_i). 
\end{equation*}
We now consider two cases: let $u \in \mathcal{B}$ satisfy either $d(u, \mathcal{M}_q(A) \cap \mathcal{B}) \leq \delta_0/4$, or $d(u, \mathcal{M}_q(A) \cap \mathcal{B}) > \delta_0/4$. Assuming the former, there exists $i \in I$ such that $\Vert v_i - u \Vert_{W^{1, 2}} \leq \delta_i/2$. In particular, the triangle inequality implies that if $\overline{v}$ is the closest element of $\mathcal{M}_q(A) \cap \mathcal{B}$ to $u$, then $\Vert \overline{v} - u \Vert_{W^{1, 2}} \leq \delta_i$. Thus, local quantitative stability implies, cf. Proposition \ref{prop: loc quant stab},  
\begin{equation*}
    \mathcal{Q}_q(u) - 1 \geq C(v_i) d_{\delta_i}(u, \mathcal{M}_q(A) \cap \mathcal{B})^{2 + \gamma_i} \geq C_0 d(u, \mathcal{M}_q(A) \cap \mathcal{B})^{2 + \gamma_0}. 
\end{equation*}
On the other hand, in the case $d(u, \mathcal{M}_q(A) \cap \mathcal{B}) \geq \delta_0/4$, we argue as follows. By (the contrapositive of) Proposition \ref{prop:precompactness extremal}, we get the existence of $\epsilon > 0$ such that
\[
    u \in \mathcal B,\quad  d(u, \mathcal{M}_q(A) \cap \mathcal{B}) >\frac{\delta_0}{4} \quad \Rightarrow \quad  \frac{1-{\rm Vol}_g(M)^{2/q-1}\|u\|_{L^2}^2}{\|\nabla u\|_{L^2}^2} \le A^{\rm opt}_q(M)-\epsilon,
\]
(notice that, as we require $d(u, \mathcal{M}_q(A) \cap \mathcal{B}) > \delta_0/4$, it follows automatically that  $u$ cannot be a constant, hence the ratio is well defined). Rearranging, the last conclusion becomes
\[
\mathcal{Q}_q(u) - 1 > \epsilon \|\nabla u\|_{L^2}^2.
\]
We claim that, there exist $l>0$ so that
\[
   u \in \mathcal B, \quad d(u, \mathcal{M}_q(A) \cap \mathcal{B}) >\frac{\delta_0}{4} \quad \Rightarrow\quad   \|\nabla u\|_{L^2}>l.
\]
If not, we could find a sequence $(u_n)\subseteq W^{1,2}(M)$ with $ \|u_n\|_{L^q}=1$ for all $n\in\N$ so that $d(u, \mathcal{M}_q(A)\cap \mathcal B) >\delta_0/4$ holds for all $n\in\N$ but $\|\nabla u_n\|_{L^2}\to 0$ as $n\uparrow\infty$. Up to a not relabelled subsequence, we thus infer that there is $u \in W^{1,2}(M)$ so that $u_n$ converges strong in $L^2$ and weak in $W^{1,2}(M)$ to $u$. However, by lower semicontinuity, it follows that $\|\nabla u\|_{L^2} =0$. In particular, $u$ must be constant by the Poincar\'e inequality, say $u\equiv c $ for some $c \in \R \setminus \{0\}$ and therefore we get directly $\|u_n\|_{W^{1,2}}\to \|u\|_{W^{1,2}}$ as $n\uparrow\infty$. Since $W^{1,2}$ is Hilbert and $u_n$ is already $W^{1,2}$-weak convergent, there is therefore strong $W^{1,2}$-convergence to the constant $c$. In particular, the Sobolev inequality implies that $u_n \rightarrow c$ in $L^{q}$ (note that this is not trivial when $q = 2^*$). In particular, $\| c\|_{L^q}=1$ and therefore $c$ is a competitor in $\mathcal{M}_q(A) \cap \mathcal{B}$. However, this would violate the uniform condition
\[
\frac{\delta_0}{4} <  d(u_n,  \mathcal{M}_q(A) \cap \mathcal{B})
 \le \frac{\|u_n - c\|_{W^{1,2}}}{\|u_n\|_{W^{1,2}}},
\]
for $n$ large enough. Finally, setting $C = \min\{C_0, l^2\epsilon\}$ we conclude the proof of Theorem \ref{thm: stability Aopt intro} (note that $d(u,\mathcal M_q(A))\le 1$ always hold).
\smallskip

\noindent Let us now discuss the straightforward modifications of the above arguments to obtain the proof of Theorem \ref{thm: stability Bopt intro}. In this case,  we simply need to consider $q=2^*,\, A={\sf S}_d^2,\, B=B_{2^*}^{\rm opt}(M)$ and $\mathcal M^1 = \mathcal M_{2^*}(B)\cap \mathcal B$. When $d(u, \mathcal{M}_{2^*}(B) \cap \mathcal{B}) \le \delta_0/4$ we repeat \emph{verbatim} the above arguments, the only difference being the use of Theorem \ref{thm:dichotomy for Bopt} and \eqref{eq:preocompact MB} replacing the role of Proposition \ref{prop:precompactness extremal}. Finally, the case $ d(u, \mathcal{M}_{2^*}(B) \cap \mathcal{B}) \geq \delta_0/4$ is now excluded by the last assumption in Theorem \ref{thm: stability Bopt intro} choosing $\delta$ accordingly.  \qed

\begin{remark}\label{rmk: no CC for Bopt}
    As already discussed in the Introduction, the concentration-compactness result in Proposition \ref{prop:precompactness extremal} made it possible to deduce in Theorem \ref{thm: stability Aopt intro} that, when $u$ is far away from the set of normalized optimizers, $\mathcal Q_q(u)$ is far away from its minimum value. This principle does not appear in \cite{DjadliDruet01} (recall Theorem \ref{thm:dichotomy for Bopt} and \eqref{eq:preocompact MB}). Consequently, the only information we get when applying our strategy in Theorem \ref{thm: stability Bopt intro} is a stability for functions that are \emph{assumed} sufficiently close to the optimizers.
\end{remark}

\section{Degenerate stability}
In this part, we prove degenerate stability results for $A^{\rm opt}_q(M)$.
\subsection{Proof of Theorem \ref{thm:degenerate Aopt main}}
Denote by $\lambda(M)$ the first non-zero eigenvalue of the Laplace-Beltrami operator on $M$, and let $\varphi \in W^{1,2}(M)$ be the corresponding eigenfunction satisfying 
\begin{equation}
-\Delta \varphi = \lambda(M) \varphi,\qquad \int \varphi \, d{\rm Vol}_g =0.
\label{eq:eigenfunction}
\end{equation}
Recall that $\lambda(M)$ admits the following variational characterization
\[
\lambda(M) = \inf_{\substack{f \in W^{1,2}(M)\\ f\neq cst.}}\frac{\| \nabla f\|_{L^2}^2}{\|f - \bar f\|_{L^2}^2},
\]
where $\bar f= \int f\, d{\rm Vol}_g$. Under the standing hypotheses, we claim that 
\begin{equation}
    \frac{q-2}{\lambda(M)}{\rm Vol}_g(M)^{\frac 2q -1} = A^{\rm opt}_q(M).
\label{eq:aoptlambda}
\end{equation}
That the left-hand side is bounded by the right-hand is follows from the well-known principle that the Sobolev inequality implies the Poincar\'e inequality on $M$, see e.g.\ the argument in \cite[Proposition 6.2.2]{BakryGentilLedoux14}. For the converse inequality, we argue as follows. By definition of $A^{\rm opt}_q(M)$, let us consider a sequence $(u_n)\subseteq W^{1,2}(M)$ of non-constant functions so that
\[
 \frac{\| u_n\|_{L^q}-{\rm Vol}_g(M)^{2/q-1}\|u_n\|_{L^2}}{\|\nabla u_n\|_{L^2}}\to A^{\rm opt}_q(M),\qquad\text{as }n\uparrow\infty.
\]
Without loss of generality, we can suppose that $\|u_n\|_{L^q}=1$ for all $n\in\N$. We are in a position to invoke Proposition  \ref{prop:precompactness extremal} to deduce that, up to a non-relabelled subsequence, $u_n \to u$ in $L^q$ and in $W^{1,2}$, for some non-zero $u\in\mathcal M_q(A)$. By assumption, we must have that $u \equiv c$ for $c ={\rm Vol}_g(M)^{-1/q}$ (as constant functions are the only extremal functions and such $c$ has unitary $L^q$-norm). By the quantitative linearization principle for Sobolev inequalities \cite[Lemma 6.7]{NobiliVIolo22} (and by suitably scaling the volume in the norms, as there it is assumed a probability reference measure), setting $f_n := u_n-c$, we deduce
\begin{align*}
    A^{\rm opt}_q(M) &= \lim_{n\uparrow\infty}  \frac{\| u_n\|_{L^q}-{\rm Vol}_g(M)^{2/q-1}\|u_n\|_{L^2}}{\|\nabla u_n\|_{L^2}} \\
    &= (q-2){\rm Vol}_g(M)^{\frac 2q-1}\lim_{n\uparrow\infty} \frac{\int |f_n - \bar f_n |^2\, d {\rm Vol}_g}{\|\nabla f_n\|_{L^2}}  \le \frac{q-2}{\lambda(M)}{\rm Vol}_g(M)^{\frac 2q-1}.
\end{align*}
This completes the proof of \eqref{eq:aoptlambda}.

Let us now consider the constant extremal $u = c = {\rm Vol}_g(M)^{- 1/q}$. By scaling, we can assume that the eigenfunction $\varphi$ in \eqref{eq:eigenfunction} satisfies $\|\varphi\|_{L^{q}}=1$. Then, the condition $\int \varphi\, d{\rm Vol}_g =0$ implies that $\varphi \in T_u \mathcal B$. Hence, we can Taylor expand the functional $\mathcal Q_q$ around $u=c$, in the direction $\varphi$ to infer 
 \begin{align*}
     \mathcal Q_q(c + \epsilon \varphi) - 1 &= \mathcal Q_q(c + \epsilon \varphi) - \mathcal Q_q(c) = \epsilon \nabla_\mathcal{B} \mathcal Q_q(c)[\varphi] + \frac{\epsilon^2}{2} \nabla^2_\mathcal{B} \mathcal Q_q(\xi)[\varphi,\varphi] \\
     &=  \epsilon^2\frac{1}{2} \nabla^2_\mathcal{B} \mathcal Q_q(c)[\varphi,\varphi] + o(1)\|\epsilon \varphi\|_{W^{1,2}}^2, 
 \end{align*}
where $\xi \in \mathcal B$ is a geodesic point between $c$ and $c+\epsilon \varphi$, and $o(1)$ goes to zero as $\epsilon \downarrow 0$. By the computation of Lemma \ref{lem: variations Q}, we know that
 \begin{align*}
     \frac{1}{2} \nabla^2_\mathcal{B} \mathcal Q_q(c)[\varphi,\varphi] &= A^{\rm opt}_q(M)\int \lambda(M) \varphi^2\, d{\rm Vol}_g + {\rm Vol}_g(M)^{\frac 2q -1}\int \varphi^2\, d{\rm Vol}_g  \\
     &\qquad -(q-1)\mathcal Q_q(c)\int c^{q-2}\varphi^2\, d{\rm Vol}_g\\
    &= (q-1){\rm Vol}_g(M)^{\frac 2q -1}\int \varphi^2\, d{\rm Vol}_g   -{\rm Vol}_g(M)^{\frac 2q -1}(q-1)\int \varphi^2\, d{\rm Vol}_g =0,
 \end{align*}
 having used \eqref{eq:aoptlambda}, and the choice $c ={\rm Vol}_g(M)^{-1/q} $. Finally, by contradiction, if $\gamma=0$ were possible, we would reach
 \[
 C \inf_{v \in \mathcal M_q(A)} \|c+\epsilon \varphi -v\|_{W^{1,2}}^2 \le  \mathcal{Q}_q(c+\epsilon \varphi)-1\approx  \|\epsilon \varphi\|^2_{W^{1,2}} \cdot o(1)
 \]
 which is false for $\epsilon \ll 1$. This concludes the proof.\qed
\subsection{Proof of Corollary \ref{cor:Frank example}}
We prove first i). Let $d>2$ and consider any $q \in (2,2^*)$. Recall  \cite{BakryEmery85,BidautVeronVeron1991,Beckner93} (also rearranging in \cite[Eq. (5)]{Frank21}) that
\[
A^{\rm opt}_q(\S^d) = \frac{q-2}{d}{\rm Vol}(\S^d)^{2/q-1},
\]
and that the only extremal functions are constant functions. Thus, Theorem \ref{thm:degenerate Aopt main} guarantees that the stability is degenerate. To check that $\gamma=2$ is optimal, we can rearrange the quantitative stability in \cite[Theorem 2]{Frank21} and we have for some $c_{d,q}>0$ 
\[
\frac{A^{\rm opt}_q(\S^d)\|\nabla u\|_{L^2}^2 - {\rm Vol}(\S^d)^{2/q-1}\| u\|_{L^2}^2 }{\|u\|_{L^q}^2 } - 1  \ge c_{d,q}  \left(\inf_{c \in \R}\frac{\| u-c\|_{W^{1,2}}}{\|u\|_{W^{1,2}}}\right)^4,
\]
for every $u \in W^{1,2}(\S^d)\setminus\{0\}$. Since in \cite[Theorem 2]{Frank21} the exponent $4$ is proven to be sharp, we directly deduce that in $\S^d$, and in the range $q \in (2,2^*)$, Theorem \ref{thm: stability Aopt intro} holds sharp with  $\gamma=2$.

We conclude by proving ii). For $d>2$, consider the closed Riemannian manifold
\begin{equation}
M_* := \mathbb S^1\left(\tfrac{1}{\sqrt{d-2}}\right)\times \mathbb S^{d-1},
\label{eq:Manifold Frank}
\end{equation}
equipped with the product metric denoted by $g$. In \cite{Frank21}, the author consider the following \textit{conformally invariant Sobolev inequality}: 
\begin{equation}
Y(M_*)\|u\|_{L^{2^*}}^2 \le \mathcal E_g(u),\qquad\forall u \in W^{1,2}(M_*)
\label{eq:conformally invariant Sobolev}
\end{equation}
where $\mathcal{E}_g(u) := \int |\nabla u|^2 + \frac{d-2}{4(d-1)}{\rm R}_g u^2\, d {\rm Vol}_g,$ and  $Y(M_*)$ is the \emph{Yamabe constant} of $M_*$. Notice
\[
{\rm R}_g = (d-2)(d-1),\qquad {\rm Vol}_g(M_*) = \frac{2\pi}{\sqrt{d-2}}{\rm Vol}(\mathbb S^{d-1}).
\]
Moreover, the Yamabe constant can be also computed explicitly (see \cite{Frank21} and \cite{Schoen89})
\[
Y(M_*) = \frac{(d-2)^2}4 {\rm Vol}_g(M_*)^{2/d}.
\]
Hence, rearranging in \eqref{eq:conformally invariant Sobolev}, we have for all $u \in W^{1,2}(M_*)$
\begin{equation}
\|u\|_{L^{2^*}}^2 \le \frac{4}{(d-2)^2}{\rm Vol}_g(M_*)^{-2/d}\|\nabla u\|_{L^2}^2 + {\rm Vol}_g(M_*)^{-2/d}\|u\|_{L^2}^2.
\label{eq:Sob Frank}
\end{equation}
Since \eqref{eq:conformally invariant Sobolev} is sharp by definition of the Yamabe constant $Y(M_*)$, we  deduce that
\[
A^{\rm opt}_{2^*}(M_*) =  \frac{4}{(d-2)^2}{\rm Vol}_g(M_*)^{-2/d}.
\]
Now, constant functions are the only extremal functions in \eqref{eq:conformally invariant Sobolev} and hence also in \eqref{eq:Sob Frank}. Thus, to apply Theorem \ref{thm:degenerate Aopt main} we need to check that
\begin{equation}
  {\sf S}_d^2 <  A^{\rm opt}_{2^*}(M_*), \label{eq:claim}   
\end{equation}
Recalling that ${\sf S}_d^2  =\frac{4}{d(d-2)}{\rm Vol}(\mathbb S^d)^{-2/d}$ and using the above computations, \eqref{eq:claim}  becomes
\[
\frac{4}{d(d-2)} < \frac{4}{(d-2)^2}(d-2)^{1/d}.
\]
which is always true for $d>2$. Thus, Theorem \ref{thm:degenerate Aopt main} guarantees that the stability is degenerate. We now check that $\gamma=2$ is optimal and conclude. By \cite[Theorem 3]{Frank21}, there is a constant $c_d>0$ so that for every $u\in W^{1,2}(M_*)\setminus\{0\}$ it holds sharp
 \begin{equation}
     \mathcal E_g(u) - Y(M_*)\|u\|_{L^{2^*}}^{2} \ge c_d \inf_{c \in\R}\frac{\mathcal E_g( u- c)^2}{\mathcal E_g( u)}. 
    \label{eq:Frank stability}
 \end{equation}
Notice that $\mathcal{E}_g(u)$ is comparable up to dimensional constants to $\|u\|_{W^{1,2}}^2$ (thanks to the fact that ${\rm R}_g$ is constant). Hence, we can rearrange  in \eqref{eq:Frank stability} using \eqref{eq:Sob Frank} to deduce
\[
\frac{A^{\rm opt}_{2^*}(M_*)\|\nabla u\|_{L^2}^2 + {\rm Vol}(M_*)^{-2/d}\|u\|_{L^2}^2}{\|u\|_{L^{2^*}}} -1 \ge c_d  \left(\inf_{c \in\R}\frac{\|u-c\|_{W^{1,2}}}{\|u\|_{W^{1,2}}}\right)^4,
\]
 for every $u\in W^{1,2}(M_*)\setminus\{0\}$, and for a possibly different constant $c_d$. Finally, this concludes the proof of ii) since the optimality in \eqref{eq:Frank stability} implies that $\gamma=2$ is sharp.\qed

\smallskip

\noindent\textbf{Acknowledgments}. F.N. is a member of INDAM-GNAMPA and acknowledges partial support of the Academy of Finland project \emph{Approximate incidence geometry}, Grant No.\ 355453 and of the MIUR Excellence Department Project awarded to the Department of Mathematics, University of Pisa, CUP I57G22000700001. F.N. would like to thank I. Y. Violo for numerous conversations around the topics of this note. Both authors would like to thank I.Y. Violo, and L. Spolaor for useful discussions and interesting comments.

\bibliographystyle{amsalpha} 

\newcommand{\etalchar}[1]{$^{#1}$}
\providecommand{\bysame}{\leavevmode\hbox to3em{\hrulefill}\thinspace}
\providecommand{\MR}{\relax\ifhmode\unskip\space\fi MR }
\providecommand{\MRhref}[2]{%
  \href{http://www.ams.org/mathscinet-getitem?mr=#1}{#2}
}
\providecommand{\href}[2]{#2}

\end{document}